\documentclass{elsart}
\usepackage{amssymb}
\usepackage{pstricks, pst-plot, pst-node, pst-tree}
\usepackage{array}
\usepackage{fancyhdr}
\usepackage{longtable}
\usepackage{tabularx}
\usepackage{multirow}
\usepackage{rotating}
\usepackage{makeidx}
\usepackage{amsmath}
\usepackage{theorem}
\usepackage[mathscr]{eucal}
\usepackage{enumerate}
\usepackage[dvips]{epsfig}
\usepackage{float}  
\journal{Applied Numerical Mathematics}
\volume{59}
\pubyear{2009}
\newcounter{mylastpage}
\issue{3--4}
\usepackage{hyperref}
\makeatletter
\def\ps@copyright{%
 \let\@oddhead\@empty
 \let\@evenhead\@empty
 \def\@oddfoot{\small\slshape\hskip-5em
   Published in \@journal\ \@volume\ (\the\@pubyear)\ no. 3--4, pp.\ \ESpagenumber{firstpage}--\ESpagenumber{mylastpage},
\href{http://dx.doi.org/10.1016/j.apnum.2008.03.012}{doi: 10.1016/j.apnum.2008.03.012}}%
 \let\@evenfoot\@oddfoot}
\makeatother
\theoremstyle{plain}
\newtheorem{Def}{Definition}[section]

\newtheorem{The}[Def]{Theorem}

\newcommand{\arsinh}{\operatorname{arsinh}}

\newcommand{\E}{\operatorname{E}}
\newcommand{\Prob}{\operatorname{P}}
\begin{document}
\begin{frontmatter}

\title{Families of efficient second order Runge-Kutta methods for the weak approximation
 of It\^{o} stochastic differential equations}

\author{Kristian Debrabant} and
\ead{debrabant@mathematik.tu-darmstadt.de}
\author{Andreas R\"{o}{\ss}ler}
\ead{roessler@mathematik.tu-darmstadt.de}
\address{Technische Universit\"{a}t Darmstadt, Fachbereich Mathematik, Schlo{\ss}gartenstr.7,
D-64289 Darmstadt, Germany}

Dedicated to Professor Karl Strehmel
\begin{abstract}
Recently, a new class of second order Runge-Kutta methods for It\^{o}
stochastic differential equations with a multidimensional Wiener
process was introduced by R\"{o}{\ss}ler \cite{Roe06c}. In contrast to
second order methods earlier proposed by other authors, this class
has the advantage that the number of function evaluations depends
only linearly on the number of Wiener processes and not
quadratically. In this paper, we give a full classification of the
coefficients of all explicit methods with minimal stage number.
Based on this classification, we calculate the coefficients of an
extension with minimized error constant of the well-known RK32
method \cite{Butcher03} to the stochastic case. For three examples, this method is
compared numerically with known order two methods and yields very
promising results.
\end{abstract}

\end{frontmatter}
\section{Introduction} \label{Introduction}
In recent years, the development of numerical methods for the
approximation of stochastic differential equations (SDEs) has become
a field of increasing interest, see. e.g \cite{KP99,Mil95} and
references therein. Whereas strong approximation methods are
designed to obtain good pathwise solutions \cite{BuBu96}, weak
approximation focuses on the expectation of functionals of the
solution. Second order stochastic Runge-Kutta (SRK) methods for the
weak approximation of SDEs were proposed by Kloeden and Platen
\cite{KP99}, Komori \cite{komori07wso}, Mackevicius and Navikas
\cite{MacNav01}, Tocino and Vigo-Aguiar \cite{ToVA02}, and the authors
 \cite{DebrabantRoessler1,Roe06b}. However,
these methods were not suitable for problems with high numbers $m$
of Wiener processes, because for these methods the number of
function evaluations per step increases quadratically in $m$.
Recently, new classes of SRK methods were introduced by R\"{o}{\ss}ler
\cite{Roe06c,Roe06d} which overcome this problem. In Section \ref{Sec:SRK-methods}
we present the one of these classes which is suitable for It\^{o} SDEs.
 The aim of this paper is to give a full
classification of all the explicit methods within this class with
minimal stage number, which is done in
Section~\ref{Sec:Parameter-families}. As an application, in Section
\ref{Sec:Optimal-schemes} we extend the well known RK32 scheme \cite{Butcher03} to
an SRK method with minimized leading local error term. The
performance of this method is illustrated by some numerical examples
in section~\ref{Sec:Numerical-Example}.
\\ \\
We denote by $(X(t))_{t \in I}$ the solution of the $d$-dimensional
It\^{o} SDE defined by
\begin{equation} \label{St-lg-sde-ito-1}
    {\mathrm{d}} X(t) = a(t,X(t)) \, {\mathrm{d}}t + b(t,X(t)) \,
    {\mathrm{d}}W(t), \qquad X({t_0}) = x_{0},
\end{equation}
with an $m$-dimensional Wiener process $(W(t))_{t \geq 0}$ and
$I=[t_0,T]$.
We assume that the Borel-measurable coefficients $a : I \times
\mathbb{R}^d \rightarrow \mathbb{R}^d$ and $b : I \times
\mathbb{R}^d \rightarrow \mathbb{R}^{d \times m}$ satisfy a
Lipschitz and a linear growth condition such that the Existence
and Uniqueness Theorem~\cite{KP99} applies. In the following, let
$b^j(t,x) = (b^{i,j}(t,x))_{1 \leq i \leq d} \in \mathbb{R}^d$
denote the $j$th column of the diffusion matrix $b(t,x)$ for $j=1,
\ldots, m$. \\ \\
Let a discretization $I_h = \{t_0, t_1, \ldots, t_N\}$ with $t_0 <
t_1 < \ldots < t_N =T$ of the time interval $I=[t_0,T]$ with step
sizes $h_n = t_{n+1}-t_n$ for $n=0,1, \ldots, N-1$ be given.
Further, define $h = \max_{0 \leq n < N} h_n$ as the maximum step
size.
Let $C_P^l(\mathbb{R}^d, \mathbb{R})$ denote the space of all $g
\in C^l(\mathbb{R}^d,\mathbb{R})$ fulfilling a polynomial growth
condition and let $g \in C_P^{k,l}(I \times \mathbb{R}^d,
\mathbb{R})$ if $g(\cdot,x) \in C^{k}(I,\mathbb{R})$ and
$g(t,\cdot) \in C_P^l(\mathbb{R}^d, \mathbb{R})$
for all $t \in I$ and $x \in \mathbb{R}^d$ \cite{KP99}.
\begin{Def}
    A time discrete approximation
    $Y=(Y(t))_{t \in I_h}$
    converges weakly
    with order $p$ to $X$ as $h \rightarrow 0$ at time $t \in I_h$ if
    for each $f \in C_P^{2(p+1)}(\mathbb{R}^d, \mathbb{R})$
    exist a constant $C_f$
    and a finite $\delta_0 > 0$ such that
    \begin{equation}
        | \E(f(X(t))) - \E(f(Y(t))) | \leq C_f \, h^p
    \end{equation}
    holds for each $h \in \, ]0,\delta_0[\,$.
\end{Def}
\section{Stochastic Runge-Kutta methods}
\label{Sec:SRK-methods}
We consider the stochastic Runge-Kutta methods introduced in
\cite{Roe06c} for the weak approximation of
SDE~(\ref{St-lg-sde-ito-1}). Therefore, the $d$-dimensional
approximation process $Y$ with $Y_n=Y(t_n)$ of an explicit $s$-stage
SRK method is defined by $Y_{0} = x_0$ and
\begin{equation} \label{SRK-method-Ito-Wm-allg01}
    \begin{split}
    Y_{n+1} = Y_n & + \sum_{i=1}^s
    \alpha_i \, a(t_n+c_i^{(0)} h_n, H_i^{(0)}) \, h_n \\
    & + \sum_{i=1}^s
    \sum_{k=1}^m
    {\beta_i^{(1)}} \, b^{k}(t_n+c_i^{(1)} h_n, H_i^{(k)}) \, \hat{I}_{(k)} \\
    & +
    \sum_{i=1}^s \sum_{k=1}^m
    {\beta_i^{(2)}} \, b^{k}(t_n+c_i^{(1)} h_n, H_i^{(k)}) \,
    \frac{\hat{I}_{(k,k)}}{\sqrt{h_n}} \\
    & +
    \sum_{i=1}^s \sum_{k=1}^m
    {\beta_i^{(3)}} \, b^{k}(t_n+{c}_i^{(2)} h_n, \hat{H}_i^{(k)}) \,
    \hat{I}_{(k)} \\
    & +
    \sum_{i=1}^s \sum_{k=1}^m
    {\beta_i^{(4)}} \, b^{k}(t_n+{c}_i^{(2)} h_n, \hat{H}_i^{(k)}) \,
    \sqrt{h_n}
    \end{split}\end{equation}
for $n=0,1, \ldots, N-1$ with stage values
\begin{alignat*}{5}
    H_i^{(0)} &&= Y_n &+ \sum_{j=1}^{i-1} A_{ij}^{(0)}
    \, a(t_n+c_j^{(0)} h_n, H_j^{(0)}) \, h_n \\
    && &+ \sum_{j=1}^{i-1} \sum_{l=1}^m
    {B_{ij}^{(0)}} \, b^l(t_n+c_j^{(1)} h_n, H_j^{(l)}) \, \hat{I}_{(l)} \\
    H_i^{(k)} &&= Y_n &+ \sum_{j=1}^{i-1} A_{ij}^{(1)}
    \, a(t_n+c_j^{(0)} h_n, H_j^{(0)}) \, h_n \\
    && &+ \sum_{j=1}^{i-1}
    {B_{ij}^{(1)}} \, b^k(t_n+c_j^{(1)} h_n, H_j^{(k)}) \,
    \sqrt{h_n} \\
    \hat{H}_i^{(k)} &&= Y_n &+ \sum_{j=1}^{s} {A}_{ij}^{(2)}
    \, a(t_n+c_j^{(0)} h_n, H_j^{(0)}) \, h_n \\
    && &+ \sum_{j=1}^{s} \sum_{\substack{l=1 \\ l \neq k}}^m
    {{B}_{ij}^{(2)}} \, b^l(t_n+c_j^{(1)} h_n, H_j^{(l)}) \,
    \frac{\hat{I}_{(k,l)}}{\sqrt{h_n}}
\end{alignat*}
for $i=1, \ldots, s$ and $k=1, \ldots, m$. Here, $\alpha,
\beta^{(1)},\dots,\beta^{(4)},c^{(q)}\in \mathbb{R}^s$ and
$A^{(q)}$, $B^{(q)} \in \mathbb{R}^{s \times s}$ for $0\leq q \leq
2$ with $A_{ij}^{(q)} = B_{ij}^{(q)} = 0$ for $j \geq i$ and $0\leq
q\leq1$ are the vectors and matrices of coefficients of the SRK
method, $c^{(q)}=A^{(q)} e$ for $0 \leq q \leq 2$ with a vector
$e=(1, \ldots, 1)^T$. In the following, the product of column
vectors is defined component-wise. The coefficients of the SRK
method~(\ref{SRK-method-Ito-Wm-allg01}) are determined by the
following Butcher tableau:
%
%
\renewcommand{\arraystretch}{1.8}
\begin{center}
\begin{tabular}{c|c|c|c}
    $c^{(0)}$ & ${A}^{(0)}$ & $B^{(0)}$ & \\
    \cline{1-4}
    $c^{(1)}$ & ${A}^{(1)}$ & $B^{(1)}$ & \\
    \cline{1-4}
    $c^{(2)}$ & ${A}^{(2)}$ & $B^{(2)}$ & \\
    \hline
    & $\alpha^T$ & ${\beta^{(1)}}^T$ & ${\beta^{(2)}}^T$ \\
    \cline{2-4}
    & & ${\beta^{(3)}}^T$ & ${\beta^{(4)}}^T$
\end{tabular}
\end{center}
\renewcommand{\arraystretch}{1.0}
$\hat{I}_{(k)}$ are three-point distributed random variables with
$\Prob(\hat{I}_{(k)} = \pm \sqrt{3 \, h_n} ) = \frac{1}{6}$ and
$\Prob(\hat{I}_{(k)} = 0 ) = \frac{2}{3}$. Further,
$\hat{I}_{(k,l)}$ are defined by
\begin{equation}
    \hat{I}_{(k,l)} = \begin{cases} \tfrac{1}{2} ( \hat{I}_{(k)}
    \hat{I}_{(l)} - \sqrt{h_n} \tilde{I}_{(k)} ) & \text{if } k <
    l \\
    \tfrac{1}{2} ( \hat{I}_{(k)}
    \hat{I}_{(l)} + \sqrt{h_n} \tilde{I}_{(l)} ) & \text{if } l < k \\
    \tfrac{1}{2} ( \hat{I}_{(k)}^2 - h_n ) & \text{if } k =
    l
    \end{cases}
\end{equation}
with two point distributed random variables $\tilde{I}_{(k)}$
satisfying $\Prob(\tilde{I}_{(k)} = \pm \sqrt{h_n})=\tfrac{1}{2}$.

By the application of the multi--colored rooted tree
analysis~\cite{Roe06a}, order conditions for the coefficients of the
SRK method~(\ref{SRK-method-Ito-Wm-allg01}) can be easily
determined.
As a result of this, the following
Theorem~\ref{SRK-theorem-ito-ord2-Wm-main1} due to
R\"{o}{\ss}ler~\cite{Roe06c} gives order conditions for the SRK
method~(\ref{SRK-method-Ito-Wm-allg01}) up to order two.
\begin{The} \label{SRK-theorem-ito-ord2-Wm-main1}
    Let $a^i, b^{i,j} \in C_P^{2,4}(I \times \mathbb{R}^d,
    \mathbb{R})$ for $1 \leq i \leq d$, $1 \leq j \leq m$.
    If the coefficients of the SRK method~(\ref{SRK-method-Ito-Wm-allg01})
    fulfill the equations
    \begin{alignat*}{5}
        1&. \quad \alpha^T e = 1 \qquad \qquad
        &2. \quad &{\beta^{(4)}}^T e = 0 \qquad \qquad
        &3. \quad &{\beta^{(3)}}^T e = 0 \\
        4&. \quad ({\beta^{(1)}}^T e)^2 = 1 \qquad \qquad
        &5. \quad &{\beta^{(2)}}^T e = 0 \qquad \qquad
        &6. \quad &{\beta^{(1)}}^T {B^{(1)}} e = 0 \\
        7&. \quad {\beta^{(4)}}^T A^{(2)} e = 0
        \qquad &8. \quad &{\beta^{(3)}}^T
        {B^{(2)}} e = 0 \qquad
        &9. \quad &{\beta^{(4)}}^T ({B^{(2)}} e)^2
        = 0
    \end{alignat*}
    then the method attains order 1 in the weak sense.
    In addition, if $a^i, b^{i,j} \in C_P^{3,6}(I \times \mathbb{R}^d,
    \mathbb{R})$ for $1 \leq i \leq d$, $1 \leq j \leq m$ and if the
    equations
    {\allowdisplaybreaks
    \begin{alignat*}{3}
        10&. \quad \alpha^T A^{(0)} e = \tfrac{1}{2}
        \qquad \qquad \qquad \qquad \quad
        &11. \quad &\alpha^T (B^{(0)} e)^2 = \tfrac{1}{2} \\
        12&. \quad ({\beta^{(1)}}^T e) (\alpha^T B^{(0)} e) =
        \tfrac{1}{2}
        \qquad \qquad
        &13. \quad &({\beta^{(1)}}^T e) ({\beta^{(1)}}^T A^{(1)} e) = \tfrac{1}{2} \\
        14&. \quad {\beta^{(3)}}^T A^{(2)} e = 0
        \qquad \qquad
        &15. \quad &{\beta^{(2)}}^T B^{(1)} e = 1 \\
        16&. \quad {\beta^{(4)}}^T B^{(2)} e = 1
        \qquad \qquad
        &17. \quad &({\beta^{(1)}}^T e) ({\beta^{(1)}}^T
        (B^{(1)} e)^2) = \tfrac{1}{2} \\
        18&. \quad ({\beta^{(1)}}^T e) ({\beta^{(3)}}^T
        (B^{(2)} e)^2) = \tfrac{1}{2} \qquad
        &19. \quad &{\beta^{(1)}}^T (B^{(1)} (B^{(1)} e)) =
        0 \\
        20&. \quad {\beta^{(3)}}^T (B^{(2)}
        (B^{(1)}e)) = 0 \qquad
        &21. \quad &{\beta^{(3)}}^T (B^{(2)} (B^{(1)}
        (B^{(1)} e))) = 0 \\
        22&. \quad {\beta^{(1)}}^T (A^{(1)} (B^{(0)} e)) =
        0 \qquad
        &23. \quad &{\beta^{(3)}}^T (A^{(2)} (B^{(0)}
        e)) = 0 \\
        24&. \quad {\beta^{(4)}}^T (A^{(2)} e)^2 = 0
        \qquad
        &25. \quad &{\beta^{(4)}}^T (A^{(2)} (A^{(0)}
        e)) = 0 \\
        26&. \quad \alpha^T (B^{(0)} (B^{(1)} e)) = 0 \qquad
        &27. \quad &{\beta^{(2)}}^T A^{(1)} e = 0 \\
        28&. \quad {\beta^{(1)}}^T ((A^{(1)} e) (B^{(1)} e)) = 0
        \qquad
        &29. \quad &{\beta^{(3)}}^T ((A^{(2)} e)
        (B^{(2)} e)) = 0 \\
        30&. \quad {\beta^{(4)}}^T (A^{(2)} (B^{(0)}
        e)) = 0 \qquad
        &31. \quad &{\beta^{(2)}}^T (A^{(1)} (B^{(0)} e))
        = 0 \\
        32&. \quad {\beta^{(4)}}^T ((B^{(2)} e)^2
        (A^{(2)} e)) = 0 \qquad
        &33. \quad &{\beta^{(4)}}^T (A^{(2)} (B^{(0)}
        e)^2) = 0 \\
        34&. \quad {\beta^{(2)}}^T (A^{(1)} (B^{(0)} e)^2) = 0
        \qquad
        &35. \quad &{\beta^{(1)}}^T (B^{(1)} ( A^{(1)} e)) =
        0 \\
        36&. \quad {\beta^{(3)}}^T (B^{(2)} (A^{(1)}
        e)) = 0 \qquad
        &37. \quad &{\beta^{(2)}}^T (B^{(1)} e)^2 = 0 \\
        38&. \quad {\beta^{(4)}}^T (B^{(2)} (B^{(1)}
        e)) = 0 \qquad
        &39. \quad &{\beta^{(2)}}^T (B^{(1)} (B^{(1)} e))
        = 0 \\
        40&. \quad {\beta^{(1)}}^T (B^{(1)} e)^3 = 0 \qquad
        &41. \quad &{\beta^{(3)}}^T (B^{(2)} e)^3 = 0 \\
        42&. \quad {\beta^{(1)}}^T (B^{(1)} (B^{(1)}
        e)^2) = 0 \qquad
        &43. \quad &{\beta^{(3)}}^T (B^{(2)} (B^{(1)}
        e)^2) = 0 \\
        44&. \quad {\beta^{(4)}}^T (B^{(2)} e)^4 = 0
        \qquad
        &45. \quad &{\beta^{(4)}}^T (B^{(2)} (B^{(1)}
        e))^2 = 0
        \\
        46&. \quad {\beta^{(4)}}^T ((B^{(2)} e)
        (B^{(2)} (B^{(1)} e))) = 0
        &47. \quad &\alpha^T ((B^{(0)} e) (B^{(0)}
        (B^{(1)} e))) = 0 \\
        48&. \quad {\beta^{(1)}}^T ((A^{(1)} (B^{(0)}
        e)) (B^{(1)} e)) = 0
        &49. \quad &{\beta^{(3)}}^T ((A^{(2)}
        (B^{(0)} e)) (B^{(2)} e)) = 0 \\
        50&. \quad {\beta^{(1)}}^T (A^{(1)} (B^{(0)} (B^{(1)} e))) = 0
        &51. \quad &{\beta^{(3)}}^T (A^{(2)} (B^{(0)}
        (B^{(1)} e))) = 0 \\
        52&. \quad {\beta^{(4)}}^T ((B^{(2)} (A^{(1)}
        e)) (B^{(2)} e)) = 0
        &53. \quad &{\beta^{(1)}}^T (B^{(1)} (A^{(1)}
        (B^{(0)} e))) = 0 \\
        54&. \quad {\beta^{(3)}}^T (B^{(2)} (A^{(1)}
        (B^{(0)} e))) = 0
        &55. \quad &{\beta^{(1)}}^T ((B^{(1)} e) (B^{(1)}
        (B^{(1)} e))) = 0 \\
        56&. \quad {\beta^{(3)}}^T ((B^{(2)} e)
        (B^{(2)} (B^{(1)} e))) = 0
        &57. \quad &{\beta^{(1)}}^T (B^{(1)} (B^{(1)}
        (B^{(1)} e))) = 0 \\
        58&. \quad {\beta^{(4)}}^T ((B^{(2)} e)
        (B^{(2)} (B^{(1)} e)^2)) = 0
        &59. \quad &{\beta^{(4)}}^T ((B^{(2)} e)
        (B^{(2)} (B^{(1)} (B^{(1)} e)))) = 0
    \end{alignat*}
}
    are fulfilled, then
    the SRK method~(\ref{SRK-method-Ito-Wm-allg01}) attains order 2
    in the weak sense.
\end{The}
It turns out that explicit order one SRK methods need at least $s=1$
stage while order two SRK methods need $s \geq 3$ stages. This is
due to e.g.\ conditions 4., 6.\ and 17., which can not be fulfilled
in the case of $s \leq 2$ stages for explicit order two SRK methods.
In the following, we distinguish between the stochastic and the
deterministic order of convergence. Let $p_S=p$ denote the order of
convergence of the SRK method if it is applied to an SDE and let
$p_D$ with $p_D \geq p_S$ denote the order of convergence of the SRK
method if it is applied to a deterministic ordinary differential
equation (ODE), i.e., SDE~(\ref{St-lg-sde-ito-1}) with $b \equiv 0$.
We also write $(p_D,p_S)$ in the following.
\section{Parameter families for SRK methods}
\label{Sec:Parameter-families}
\subsection{Coefficients for SRK methods of order (1,1)}
\label{Sec:Coeff-SRK-Ord-1-1}
First, we analyze explicit SRK
methods~(\ref{SRK-method-Ito-Wm-allg01}) of order $p_D=p_S=1$ with
$s=1$ stage. Considering the order one conditions 1.--9.\ in
Theorem~\ref{SRK-theorem-ito-ord2-Wm-main1}, the corresponding
coefficients are uniquely determined for $c_1 \in \{-1,1\}$ by
\begin{equation} \label{Parameter-Ord11-all}
    \alpha_1 = 1, \qquad \beta_1^{(1)} = c_1,
    \qquad \beta_1^{(2)} = 0, \qquad
    \beta_1^{(3)} = 0, \qquad
    \beta_1^{(4)} = 0.
\end{equation}
The resulting class of SRK schemes coincides with the well-known
Euler-Maruyama scheme.
\subsection{Coefficients for SRK methods of order (2,1)}
\label{Sec:Coeff-SRK-Ord-2-1}
Next, we consider the case of $s=2$ stage explicit SRK methods
(\ref{SRK-method-Ito-Wm-allg01}). As already mentioned in
Section~\ref{Sec:SRK-methods}, it is not possible to attain order
$p_S=2$. However, we can find some SRK methods of order $p_D=2$ and
$p_S=1$ corresponding to the following parameter family: From
condition 1.\ of Theorem~\ref{SRK-theorem-ito-ord2-Wm-main1} follows
$\alpha_1 = 1 -\alpha_2$ and taking into account the order two
condition 10.\ we obtain $\alpha_2 = \frac{1}{2 A^{(0)}_{21}}$ for
$A^{(0)}_{21} \neq 0$. Further, condition 2.\ yields $\beta^{(4)}_1
= -\beta^{(4)}_2$, condition 3.\ results in $\beta^{(3)}_1 =
-\beta^{(3)}_2$ and condition 5.\ is fulfilled if $\beta^{(2)}_1 =
-\beta^{(2)}_2$ while condition 4.\ holds for $\beta^{(1)}_1 = c_1 -
\beta^{(1)}_2$ with $c_1 \in \{-1,1\}$. Finally, considering
condition 6.\ we need that $\beta^{(1)}_2 = 0$ or $B^{(1)}_{21}=0$,
considering condition 8.\ analogously that $\beta^{(3)}_2 = 0$ or
$B^{(2)}_{11}+B^{(2)}_{12}=B^{(2)}_{21}+B^{(2)}_{22}$  and for
condition 7.\ and 9.\ that $\beta^{(4)}_2 = 0$ or
$A^{(2)}_{11}+A^{(2)}_{12}=A^{(2)}_{21}+A^{(2)}_{22}$ and
$(B^{(2)}_{11}+B^{(2)}_{12})^2=(B^{(2)}_{21}+B^{(2)}_{22})^2$ hold.
Thus, this class of SRK methods is determined by
\begin{alignat}{7}
    \alpha^T &= \begin{bmatrix} 1-\frac{1}{2 c_2} && \frac{1}{2
    c_2} \end{bmatrix} , &\quad \quad
    {\beta^{(1)}}^T &= \begin{bmatrix} c_1-c_4 && c_4 \end{bmatrix},
    &\quad \quad
    {\beta^{(2)}}^T &= \begin{bmatrix} c_5 && -c_5 \end{bmatrix}, \notag \\
    && {\beta^{(3)}}^T &= \begin{bmatrix} c_6 && -c_6 \end{bmatrix},
    &\quad
    {\beta^{(4)}}^T &= \begin{bmatrix} c_7 && -c_7 \end{bmatrix}, \notag \\
    A^{(0)} &= \begin{bmatrix} 0 && 0 \\ c_2 && 0 \end{bmatrix}, &\quad
    A^{(1)} &= \begin{bmatrix} 0 && 0 \\ c_8 && 0 \end{bmatrix}, &\quad
    A^{(2)} &= \begin{bmatrix} c_9 && c_{10} \\ c_{11} && c_{12} \end{bmatrix},
    \notag \\
    B^{(0)} &= \begin{bmatrix} 0 && 0 \\ c_3 && 0 \end{bmatrix},
    &\quad
    B^{(1)} &= \begin{bmatrix} 0 && 0 \\ c_{17} && 0 \end{bmatrix}, &\quad
    B^{(2)} &= \begin{bmatrix} c_{13} && c_{14} \\ c_{15} && c_{16} \label{Parameter-Ord21-all}
    \end{bmatrix},
\end{alignat}
for $c_1 \in \{-1,1\}$ and $c_2, \ldots, c_{17} \in \mathbb{R}$ with
$c_2 \neq 0$, $c_4 \, c_{17} =0$, $c_6(c_{13}+c_{14}-c_{15}-c_{16})
= 0$, $c_7(c_{9}+c_{10}-c_{11}-c_{12}) = 0$ and
$c_7((c_{13}+c_{14})^2-(c_{15}+c_{16})^2) = 0$.
\subsection{Coefficients for SRK methods of order (2,2)}
\label{Sec:Coeff-SRK-Ord-2-2}
Now, we consider explicit SRK
methods~(\ref{SRK-method-Ito-Wm-allg01}) of order $p_D=p_S=2$ with
$s=3$ stages. Then, the SRK schemes of the class under consideration
are completely characterized by the following families of
coefficients which follow from the order conditions in
Theorem~\ref{SRK-theorem-ito-ord2-Wm-main1}:

We have $\alpha_1 = 1 - \alpha_2 - \alpha_3$ due to condition 1.
From condition 4. it follows that
$\beta^{(1)}_3=c_1-\beta^{(1)}_2-\beta^{(1)}_1$ with
$c_1\in\{-1,1\}$. Due to conditions 2., 3., 7., 24., 16., 14., 18.
and 8. we need
$\sum_{i=1}^3A^{(2)}_{1i}=\sum_{i=1}^3A^{(2)}_{2i}=\sum_{i=1}^3A^{(2)}_{3i}$.
From conditions 3., 8., 18. and 41. follows that
$\sum_{i,j=1}^3B^{(2)}_{i,j}=0$
 and that ${{b}}_i:=\sum_{j=1}^3B^{(2)}_{ij}$, $i=1,2,3$ are pairwise different.
Further, we have
$\beta^{(3)}_1=\frac{c_1}{2({{b}}_1-{{b}}_2)({{b}}_2+2{{b}}_1)}$,
$\beta^{(3)}_2=\frac{c_1}{2({{b}}_2-{{b}}_1)({{b}}_1+2{{b}}_2)}$,
$\beta^{(3)}_3=\frac{c_1}{2(2{{b}}_1+{{b}}_2)({{b}}_1+2{{b}}_2)}$.
From conditions 2., 9. and 16. we obtain
$\beta^{(4)}_1=\frac{{{b}}_1}{({{b}}_1-{{b}}_2)({{b}}_2+2{{b}}_1)}$,
$\beta^{(4)}_2=\frac{{{b}}_2}{({{b}}_2-{{b}}_1)({{b}}_1+2{{b}}_2)}$,
$\beta^{(4)}_3=\frac{-{{b}}_1-{{b}}_2}{(2{{b}}_1+{{b}}_2)({{b}}_1+2{{b}}_2)}$.
With 44. it follows now from the above that exactly for one $i$ from
$1,2,3$ it holds ${{b}}_i=0$. Without loss of generality we can
assume that ${{b}}_1=0$. Due to conditions 15. and 37. we need
$B^{(1)}_{21} \neq 0$ and from conditions 17. and 40. follows
$\beta^{(1)}_3 \neq 0$. Now, by condition 19. follows that
$B^{(1)}_{32}=0$ and we deduce from 6., 17. and 40. that
$B^{(1)}_{31} = -B^{(1)}_{21}$. With conditions 5., 15. and 37.
follows that $\beta^{(2)}_2=\frac1{2B^{(1)}_{21}}$,
$\beta^{(2)}_3=-\frac1{2B^{(1)}_{21}}$ and $\beta^{(2)}_1=0$.
Conditions 4., 6. and 17. yield
$\beta^{(1)}_1=c_1-\frac{c_1}{2(B^{(1)}_{21})^2}$ and
$\beta^{(1)}_2=\beta^{(1)}_3=\frac{c_1}{4(B^{(1)}_{21})^2}$. From
condition 56. and 58. we obtain $B^{(2)}_{22}=-B^{(2)}_{33}$ and
$B^{(2)}_{23}=-B^{(2)}_{33}$. Condition 43. yields now
$B^{(2)}_{12}=-B^{(2)}_{13}$. Condition 46. gives
$B^{(2)}_{32}=B^{(2)}_{33}$. Condition 20. leads to
$B^{(2)}_{13}=0$. Now, due to 13. and 28. we need that $A^{(1)}_{21}
= (B^{(1)}_{21})^2$ and $A^{(1)}_{31} = (B^{(1)}_{21})^2 -
A^{(1)}_{32}$.
\\
To fulfill 11., 12., 22., 23., 30. and 33., we obtain the following cases:
\begin{enumerate}[1)]
\item\label{PF1} $B^{(0)}_{21}=0,B^{(0)}_{31}+B^{(0)}_{32}=c_1,\alpha_3=\frac12,A^{(2)}_{13}=A^{(2)}_{23}=A^{(2)}_{33}$,
\item\label{PF2}
$A^{(1)}_{32}=0,B^{(0)}_{21}=c_1,B^{(0)}_{31}+B^{(0)}_{32}=0,\alpha_2=\frac12,A^{(2)}_{12}=A^{(2)}_{22}=A^{(2)}_{32}$,
\item\label{PF3}
$A^{(1)}_{32}=0,B^{(0)}_{21}=B^{(0)}_{31}+B^{(0)}_{32}=c_1,\alpha_2+\alpha_3=\frac12,A^{(2)}_{12}+A^{(2)}_{13}=A^{(2)}_{22}+A^{(2)}_{23}=A^{(2)}_{32}+A^{(2)}_{33}$,
\item\label{PF4} $A^{(1)}_{32}=0,B^{(0)}_{21}\neq0\neq B^{(0)}_{31}+B^{(0)}_{32}\neq B^{(0)}_{21},
A^{(2)}_{22}=A^{(2)}_{32},A^{(2)}_{23}=A^{(2)}_{33},A^{(2)}_{12}=A^{(2)}_{32}+(A^{(2)}_{33}-A^{(2)}_{13})\frac{B^{(0)}_{31}+B^{(0)}_{32}}{B^{(0)}_{21}}$,
$\alpha_2=\frac12\frac{1-c_1(B^{(0)}_{31}+B^{(0)}_{32})}{B^{(0)}_{21}(B^{(0)}_{21}-B^{(0)}_{31}-B^{(0)}_{32})}$,
$\alpha_3=-\frac12\frac{1-c_1B^{(0)}_{21}}{(B^{(0)}_{31}+B^{(0)}_{32})(B^{(0)}_{21}-B^{(0)}_{31}-B^{(0)}_{32})}$.
\end{enumerate}
However, from 26. and 51. it follows
\begin{enumerate}[a)]
\item\label{PFa} $B^{(0)}_{32}=0$ or
\item\label{PFb} $\alpha_3=0,A^{(2)}_{23}+A^{(2)}_{33}=2A^{(2)}_{13}$.
\end{enumerate}
Finally, the equations 10. and 25. imply the cases
\begin{enumerate}[i)]
\item\label{PFi}
$\alpha_2(A^{(2)}_{23}-A^{(2)}_{33})\neq\alpha_3(A^{(2)}_{22}-A^{(2)}_{32})$,
$A^{(0)}_{21}=\frac12\frac{A^{(2)}_{23}-A^{(2)}_{33}}{\alpha_2(A^{(2)}_{23}-A^{(2)}_{33})-\alpha_3(A^{(2)}_{22}-A^{(2)}_{32})}$,
$A^{(0)}_{31}=-\frac12\frac{A^{(2)}_{22}-A^{(2)}_{32}}{\alpha_2(A^{(2)}_{23}-A^{(2)}_{33})-\alpha_3(A^{(2)}_{22}-A^{(2)}_{32})}-A^{(0)}_{32}$,
\item\label{PFii}
$A^{(2)}_{23}=A^{(2)}_{33}$, $A^{(2)}_{22}=A^{(2)}_{32}$, $\alpha_2\neq0$, $A^{(0)}_{21}=\frac{1-2\alpha_3(A^{(0)}_{31}+A^{(0)}_{32})}{2\alpha_2}$,
\item\label{PFiii}
$A^{(2)}_{23}=A^{(2)}_{33}$, $A^{(2)}_{22}=A^{(2)}_{32}$, $\alpha_2=0$, $\alpha_3\neq0$, $A^{(0)}_{31}=\frac1{2\alpha_3}-A^{(0)}_{32}$.
\end{enumerate}
With these settings, all the remaining order conditions are now fulfilled.\\
Summarizing our results, we have the
following classification for the SRK schemes of order $p_D=p_S=2$
for the considered class with $s=3$ stages: For $c_1 \in \{-1,1\}$
and $c_2,c_3,c_4,c_5 \in \mathbb{R}$ with $c_3 \neq 0$ and $c_4 \neq
0$ holds
\begin{alignat}{5}
    {\beta^{(1)}}^T &= \begin{bmatrix} c_1-\frac{c_1}{2 c_3^2} & &
    \frac{c_1}{4 c_3^2} & & \frac{c_1}{4 c_3^2} \end{bmatrix},
    &\quad \quad
    {\beta^{(2)}}^T &= \begin{bmatrix} 0 & & \frac{1}{2 c_3} & &
    -\frac{1}{2 c_3} \end{bmatrix}, \\
    {\beta^{(4)}}^T &= \begin{bmatrix} 0 & & \frac{1}{2 c_4} & &
    -\frac{1}{2 c_4} \end{bmatrix},
    &\quad \quad
        {\beta^{(3)}}^T &= \begin{bmatrix} -\frac{c_1}{2 c_4^2} & &
    \frac{c_1}{4 c_4^2} & & \frac{c_1}{4 c_4^2} \end{bmatrix},\\
    A^{(1)} &= \begin{bmatrix} 0 && 0 && 0 \\ c_3^2 && 0 && 0 \\
    c_3^2-c_2 && c_2 && 0 \end{bmatrix}, &\quad \quad B^{(1)} &=
    \begin{bmatrix} 0 && 0 && 0 \\ c_3 && 0 && 0 \\ -c_3 && 0 && 0
    \end{bmatrix}, \label{Parameter-Ord22-A1-B1} \\
    B^{(2)} &=
    \begin{bmatrix} 0 && 0 && 0 \\ c_4+2c_5 && -c_5 && -c_5 \\ -c_4 -2c_5&& c_5 && c_5
    \end{bmatrix} . \label{Parameter-Ord22-A2-B2}
\end{alignat}
Now, the following cases are possible: \\ \\
In the case \ref{PF1}\ref{PFa}\ref{PFi}) we get with
$c_6,\dots,c_{12}\in\mathbb{R}$ that
\begin{equation} 
    \alpha^T = \begin{bmatrix} \frac12-c_{11} & &
    c_{11} & & \frac12 \end{bmatrix} ,
    \quad
    B^{(0)} =
    \begin{bmatrix} 0 && 0 && 0 \\ 0 && 0 && 0 \\c_1  &&0  && 0
    \end{bmatrix},
\end{equation}
\begin{equation}
    A^{(0)} = \begin{bmatrix} 0 && 0 && 0 \\ 0 && 0 && 0 \\
    c_{12} && 1-c_{12}&& 0 \end{bmatrix},
\quad
A^{(2)} = \begin{bmatrix} c_6-c_7 && c_7 && c_8 \\ c_6-c_9 && c_9 && c_{8} \\
    c_6-c_{10}&& c_{10} && c_{8} \end{bmatrix}.
\end{equation}
In the case \ref{PF1}\ref{PFa}\ref{PFii}) we obtain with
$c_6,\dots,c_{12}\in\mathbb{R}$ and $c_{10}\neq0$ that
\begin{equation} 
    \alpha^T = \begin{bmatrix} \frac12-c_{10} & &
    c_{10} & & \frac12 \end{bmatrix} ,
    \quad
    B^{(0)} =
    \begin{bmatrix} 0 && 0 && 0 \\ 0 && 0 && 0 \\c_1  &&0  && 0
    \end{bmatrix},
\end{equation}
\begin{equation}
    A^{(0)} = \begin{bmatrix} 0 && 0 && 0 \\\frac{1-c_{11}}{2c_{10}}  && 0 && 0 \\
    c_{11}-c_{12} && c_{12}&& 0 \end{bmatrix},
\quad
A^{(2)} = \begin{bmatrix} c_6-c_7 && c_7 && c_8 \\ c_6-c_9 && c_9 && c_{8} \\
    c_6-c_{9}&& c_{9} && c_{8} \end{bmatrix}.
\end{equation}
Considering the case \ref{PF1}\ref{PFa}\ref{PFiii}) we obtain with
$c_6,\dots,c_{11}\in\mathbb{R}$ that
\begin{equation} 
    \alpha^T = \begin{bmatrix} \frac12 & &
    0 & & \frac12 \end{bmatrix} ,
    \quad
    B^{(0)} =
    \begin{bmatrix} 0 && 0 && 0 \\ 0 && 0 && 0 \\c_1  &&0  && 0
    \end{bmatrix},
\end{equation}
\begin{equation}
    A^{(0)} = \begin{bmatrix} 0 && 0 && 0 \\ c_{10} && 0 && 0 \\
    1-c_{11} && c_{11}&& 0 \end{bmatrix},
\quad
A^{(2)} = \begin{bmatrix} c_6-c_7 && c_7 && c_8 \\ c_6-c_9 && c_9 && c_{8} \\
    c_6-c_{9}&& c_{9} && c_{8} \end{bmatrix}.
\end{equation}
For the case \ref{PF2}\ref{PFa}\ref{PFi}) we get with $c_2=0$ in
(\ref{Parameter-Ord22-A1-B1}) and $c_6,\dots,c_{12}\in \mathbb{R}$,
$c_9\neq c_{10}$ the coefficients
\begin{equation} 
    \alpha^T = \begin{bmatrix} \frac{1}{2}-c_{11} &&
    \frac12 && c_{11} \end{bmatrix},
    \quad
    B^{(0)} =
    \begin{bmatrix} 0 && 0 && 0 \\ c_1 && 0 && 0 \\ 0 && 0 && 0
    \end{bmatrix},
\end{equation}
\begin{equation}
    A^{(0)} = \begin{bmatrix} 0 && 0 && 0 \\ 1 && 0 && 0 \\
    c_{12} && -c_{12}&& 0 \end{bmatrix},
\quad
A^{(2)} = \begin{bmatrix} c_6-c_8 && c_7 && c_8 \\ c_6-c_{9} && c_7 && c_{9} \\
    c_6-c_{10} && c_{7} && c_{10} \end{bmatrix}.
\end{equation}
In the case \ref{PF2}\ref{PFa}\ref{PFii}) we obtain with $c_2=0$ in
(\ref{Parameter-Ord22-A1-B1}) and $c_6,\dots,c_{12}\in \mathbb{R}$
the coefficients
\begin{equation}
    \alpha^T = \begin{bmatrix} \frac{1}{2}-c_{10} &&
    \frac12 && c_{10} \end{bmatrix},
    \quad
    B^{(0)} =
    \begin{bmatrix} 0 && 0 && 0 \\ c_1 && 0 && 0 \\ 0 && 0 && 0
    \end{bmatrix},
\end{equation}
\begin{equation}
    A^{(0)} = \begin{bmatrix} 0 && 0 && 0 \\ 1-2c_{10}c_{11} && 0 && 0 \\
    c_{11}-c_{12} && c_{12}&& 0 \end{bmatrix},
\quad
A^{(2)} = \begin{bmatrix} c_6-c_8 && c_7 && c_8 \\ c_6-c_{9} && c_7 && c_{9} \\
    c_6-c_{9} && c_{7} && c_{9} \end{bmatrix}.
\end{equation}
For the case \ref{PF2}\ref{PFb}\ref{PFi}) we get with $c_2=0$ in
(\ref{Parameter-Ord22-A1-B1}) and $c_6,\dots,c_{11}\in \mathbb{R}$,
$c_8\neq c_9$ the coefficients
\begin{equation}
    \alpha^T = \begin{bmatrix} \frac{1}{2}&&
    \frac12 && 0\end{bmatrix},
    \quad B^{(0)} =
    \begin{bmatrix} 0 && 0 && 0 \\ c_1 && 0 && 0 \\ c_{10} && -c_{10} && 0
    \end{bmatrix},
\end{equation}
\begin{equation}
    A^{(0)} = \begin{bmatrix} 0 && 0 && 0 \\
    1 && 0 && 0 \\
    c_{11}&&-c_{11}&& 0 \end{bmatrix},
\quad
A^{(2)} = \begin{bmatrix} c_6-c_8 && c_7 && c_8 \\ c_6-c_{9} && c_7 && c_{9} \\
    c_6-2c_8+c_9 && c_{7} && 2c_8-c_9 \end{bmatrix}.
\end{equation}
In the cases \ref{PF2}\ref{PFb}\ref{PFii}),
\ref{PF3}\ref{PFb}\ref{PFii}) respectively \ref{PF4}\ref{PFb}\ref{PFii})
we obtain with $c_2=0$ in (\ref{Parameter-Ord22-A1-B1}) and
 $c_6,\dots, c_{12} \in\mathbb{R}$
\begin{equation}
    \alpha^T = \begin{bmatrix}
    \frac12
    && \frac12
    &&0
        \end{bmatrix}, \,\,
\,\, B^{(0)} =
    \begin{bmatrix} 0 && 0 && 0 \\
    c_{1}&& 0 && 0 \\
    c_{9}-c_{10} && c_{10} && 0
    \end{bmatrix} ,
\end{equation}
\begin{equation}
   A^{(0)} = \begin{bmatrix} 0 && 0 && 0 \\ 1 && 0 && 0 \\
    c_{11} && c_{12} && 0 \end{bmatrix},
\quad
A^{(2)} = \begin{bmatrix} c_6 && c_7 && c_8
\\ c_6 && c_7 && c_{8} \\
    c_6&& c_{7} && c_{8} \end{bmatrix},
\end{equation}
where $c_9=0$ in the case \ref{PF2}\ref{PFb}\ref{PFii}), $c_9=c_1$ for \ref{PF3}\ref{PFb}\ref{PFii}) and $0\neq c_9\neq c_1$ for \ref{PF4}\ref{PFb}\ref{PFii}), respectively.\\
Considering the case \ref{PF3}\ref{PFa}\ref{PFi}) we find
with $c_2=0$ in
(\ref{Parameter-Ord22-A1-B1}) and $c_6,\dots,c_{12}\in\mathbb{R}$, $c_9\neq c_{10}$ that
\begin{equation} 
    \alpha^T = \begin{bmatrix} \frac{1}{2}
    && \frac12-c_{11} && c_{11} \end{bmatrix},
    B^{(0)} =
    \begin{bmatrix} 0 && 0 && 0 \\ c_1 && 0 && 0 \\ c_1 && 0 && 0
    \end{bmatrix} ,
\end{equation}
\begin{equation}
    A^{(0)} = \begin{bmatrix} 0 && 0 && 0 \\ 1 && 0 && 0 \\
    1-c_{12} && c_{12} && 0 \end{bmatrix},
    \,\, \quad
A^{(2)} = \begin{bmatrix} c_6&& c_7 && c_8 \\ c_6&& c_9 && c_7+c_8-c_{9} \\
    c_6&& c_{10} && c_7+c_8-c_{10} \end{bmatrix}.
\end{equation}
For the case \ref{PF3}\ref{PFa}\ref{PFii}) we obtain with $c_2=0$ in
(\ref{Parameter-Ord22-A1-B1}) and $c_6,\dots,c_{12}\in\mathbb{R}$,
$c_{10}\neq\frac12$, that
\begin{equation} 
    \alpha^T = \begin{bmatrix} \frac{1}{2}
    && \frac12-c_{10} && c_{10} \end{bmatrix},
    B^{(0)} =
    \begin{bmatrix} 0 && 0 && 0 \\ c_1 && 0 && 0 \\ c_1 && 0 && 0
    \end{bmatrix} ,
\end{equation}
\begin{equation}
    A^{(0)} = \begin{bmatrix} 0 && 0 && 0 \\ \frac{1-2c_{10}c_{11}}{1-2c_{10}} && 0 && 0 \\
    c_{11}-c_{12} && c_{12} && 0 \end{bmatrix},
    \,\, \quad
A^{(2)} = \begin{bmatrix} c_6&& c_7 && c_8 \\ c_6&& c_9 && c_7+c_8-c_{9} \\
    c_6&& c_{9} && c_7+c_8-c_{9} \end{bmatrix}.
\end{equation}
In the case \ref{PF3}\ref{PFa}\ref{PFiii}) respectively \ref{PF4}\ref{PFa}\ref{PFiii}) we get with $c_2=0$ in (\ref{Parameter-Ord22-A1-B1}) and
 $c_6,\dots, c_{12} \in\mathbb{R}$ with $0\neq c_{10}$
\begin{equation}
    \alpha^T = \begin{bmatrix}
    \frac12
    &&0
    &&\frac12
        \end{bmatrix}, \,\,
\,\, B^{(0)} =
    \begin{bmatrix} 0 && 0 && 0 \\
    c_{10}&& 0 && 0 \\
    c_{1}&&0 && 0
    \end{bmatrix} ,
\end{equation}
\begin{equation}
    A^{(0)} = \begin{bmatrix} 0 && 0 && 0 \\ c_{11} && 0 && 0 \\
    1-c_{12} && c_{12} && 0 \end{bmatrix},
\,\,
A^{(2)} = \begin{bmatrix} c_6+(c_9-c_7)(1-\frac{c_{1}}{c_{10}}) && c_8+(c_{9}-c_7)\frac{c_{1}}{c_{10}} && c_7
\\ c_6&& c_8 && c_{9} \\
    c_6&& c_{8} && c_{9} \end{bmatrix},
\end{equation}
where $c_{10}=c_1$ in the case \ref{PF3}\ref{PFa}\ref{PFiii}) and $c_{10}\neq c_1$ in the case \ref{PF4}\ref{PFa}\ref{PFiii}).\\
Considering the case \ref{PF3}\ref{PFb}\ref{PFi}) we obtain
with $c_2=0$ in
(\ref{Parameter-Ord22-A1-B1}) and $c_6,\dots,c_{11}\in\mathbb{R}$, $c_7\neq c_9$, that
\begin{equation}
    \alpha^T = \begin{bmatrix} \frac{1}{2}
    && \frac12 && 0 \end{bmatrix},
    B^{(0)} =
    \begin{bmatrix} 0 && 0 && 0 \\ c_1 && 0 && 0 \\ c_1-c_{10} && c_{10} && 0
    \end{bmatrix} ,
\end{equation}
\begin{equation}
    A^{(0)} = \begin{bmatrix} 0 && 0 && 0 \\ 1 && 0 && 0 \\
    1-c_{11} && c_{11} && 0 \end{bmatrix},
    \,\, \quad
A^{(2)} = \begin{bmatrix} c_6&& c_7 && c_8 \\ c_6&& c_9 && c_7+c_8-c_{9} \\
    c_6&& 2c_7-c_9 && c_9+c_8-c_{7} \end{bmatrix}.
\end{equation}
\\
Next, we have the case \ref{PF4}\ref{PFa}\ref{PFii}) with $c_2=0$
in (\ref{Parameter-Ord22-A1-B1}) and
 $c_6,\dots, c_{13} \in\mathbb{R}$ with $0\neq c_{10}\neq c_{11}\neq0$, $c_{11}\neq c_1$,
\begin{subequations}\label{PF4aii}
\begin{equation}
    \alpha^T = \begin{bmatrix}
    1+\frac{1-c_1(c_{10}+c_{11})}{2c_{10}c_{11}}
    && \frac12\frac{1-c_1 c_{11}}{c_{10}(c_{10}-c_{11})}
    &&-\frac12\frac{1-c_1 c_{10}}{c_{11}(c_{10}-c_{11})}
        \end{bmatrix},
\end{equation}
\begin{equation}
 B^{(0)} =
    \begin{bmatrix} 0 && 0 && 0 \\
    c_{10}&& 0 && 0 \\
    c_{11}&&0 && 0
    \end{bmatrix} ,\quad
    A^{(0)} = \begin{bmatrix} 0 && 0 && 0 \\ \frac{c_{10}}{c_{11}}\frac{c_{11}(c_{11}-c_{10})-c_{12}(1-c_1c_{10})}{c_1c_{11}-1} && 0 && 0 \\
    c_{12}-c_{13} && c_{13} && 0 \end{bmatrix}
\end{equation}
\begin{equation}
A^{(2)} = \begin{bmatrix} c_6+(c_9-c_7)(1-\frac{c_{11}}{c_{10}}) && c_8+(c_{9}-c_7)\frac{c_{11}}{c_{10}} && c_7
\\ c_6&& c_8 && c_{9} \\
    c_6&& c_{8} && c_{9} \end{bmatrix}.
\end{equation}
\end{subequations}
In the remaining cases \ref{PF1}\ref{PFb}\ref{PFi})-\ref{PF1}\ref{PFb}\ref{PFiii}), \ref{PF2}\ref{PFa}\ref{PFiii}),
\ref{PF2}\ref{PFb}\ref{PFiii}),\ref{PF3}\ref{PFb}\ref{PFiii}), \ref{PF4}\ref{PFa}\ref{PFi}),
\ref{PF4}\ref{PFb}\ref{PFi}) and \ref{PF4}\ref{PFb}\ref{PFiii}) there doesn't exist a solution.
\section{Application: An SRK scheme with minimized error coefficients}
\label{Sec:Optimal-schemes}
Based on the classification given in section \ref{Sec:Parameter-families},
as an example we will now extend the well
known method RK32 of Kutta \cite{Butcher03} to an
SRK method of order (3,2).  The Butcher array of RK32
is obtained from family (\ref{PF4aii}) by setting $c_{10}=\frac{6\mp\sqrt{6}}{10}$,
$c_{11}=\frac{3\pm2\sqrt{6}}{5}$, $c_{12}=1$, $c_{13}=2$.
Due to some symmetry in the method, the sign of $c_1$ has
no influence, and we choose $c_1=1$.
Now, we want to determine the remaining coefficients by
minimizing the expectation of the local error.
Therefore, we distinguish between the cases
$m=1$ (only one Wiener process) and $m>1$. In the case of
$m=1$, to save computational effort we require that $A^{(2)}$
equals the zero matrix, because then we don't have to evaluate
the stages $\hat{H}^{(k)}$.
\\
Now, let $le_f(h)$ be the weak local error of the method starting
at the point $(t,x)$
with respect to the functional $f$ and step size $h$,
i.\,e.
\[
le_f(h)= \E \big(f(Y(t+h))-f(X(t+h))|Y(t)=X(t)=x\big).
\]
As in the deterministic case, by the colored rooted tree analysis one obtains
the representation
\[
le_f(h)=\sum_{\substack{ \textbf{t} \in
TS(\Delta) \\ \rho(\textbf{t})=3}} lec_{\textbf{t}} \,
F(\textbf{t})(x) \, h^3+\mathcal{O}(h^4),
\]
where $TS(\Delta)$ denotes a set of trees, $\rho(\textbf{t})$
the order of the tree $\textbf{t}$, $F(\textbf{t})$ the elementary differential
connected with the tree $\textbf{t}$ and $lec_{\textbf{t}}$ a coefficient
depending only on $\textbf{t}$ and the numerical method
(see \cite{Roe06a,Roe06b,Roe06c} for details).\\
Let $lec=(lec_{\textbf{t}})_{\textbf{t} \in
TS(\Delta)}$ be the vector of these coefficients. In the following, we want to
minimize $\|lec\|$ in the Euclidean norm.
Then, using again the rooted tree analysis, a tedious calculation (for
$m=1$ there exist 164 rooted trees of order three) yields that in
the Euclidean norm we have in the case of
$c_{10}=\frac{6-\sqrt{6}}{10}$
\\\centerline{$
\|lec\|^2=
\frac{60500644673 + 24530366872\sqrt6 - (217 + 88\sqrt6)(128250000c_3^2
- 92062500c_3^4)}{24000000(24 + 11\sqrt6)^2}
$}\\
which is minimized by $c_3=\pm3\sqrt\frac{38}{491}$, which gives
$\|lec\|\approx 1.275$. In the case of
$c_{10}=\frac{6+\sqrt{6}}{10}$ instead, we would obtain
$\|lec\|\approx 1.296$, so we choose the minus sign in the following.
The remaining coefficients of the method, $c_4$ and $c_5$, are determined by considering $\|lec\|$ in the case of two Wiener processes. Its
minimal value $2.859$ is attained for $c_5=\pm \frac{491}{513}\sqrt\frac{221}{4955}$,
$c_4=\mp4\sqrt\frac{221}{4955}$. For $c_3$, $c_4$ and $c_5$, the method
is invariant to the choice of the sign, so we obtain finally the scheme
DRI1 presented in Table \ref{Table-Coeff-DRI1}.

In the case $m>1$, we cannot avoid completely the evaluation of the stages $\hat{H}^{(k)}$ by letting $A^{(2)}$ equal to the zero matrix, so one could try to use the
additional degrees of freedom to minimize the local error. The
resulting method differs from DRI1 only in $A^{(2)}$, which is now
given by
\[
A^{(2)}=
\left(
\begin{matrix}
\frac{2(-3442595658 + 1259007085\sqrt6)}{1554073317(-6+\sqrt6)} & &
-\frac{8(212963260 +73915807\sqrt6)}{1554073317(-6\sqrt6)} &&
\frac{(4(-1111473969 + 371403611\sqrt6)}{23311099755}\\
\frac2{27}(7 - 2\sqrt6) & &\frac8{81}(3+\sqrt6)&&\frac4{81}(-3+\sqrt6) \\
\frac2{27}(7 - 2\sqrt6) & &\frac8{81}(3+\sqrt6) & &\frac4{81}(-3+\sqrt6)
\end{matrix}\right),
\]
and by $c^{(2)}$ which fulfills now
$c^{(2)}=(\frac23,\frac23,\frac23)$. For $m=1$, this method has
again $\|lec\|\approx 1.275$. But in the case of $m=2$, we
achieve $\|lec\|\approx 2.765$. However, this is only $3.4\%$
less than $\|lec\|$ achieved by DRI1. Due to the additional $m$
function evaluations needed for $A^{(2)}\neq0$ in the case $m>1$ (because for $A^{(2)}=0$ we would have $\hat{H}_1^{(k)}=H_1^{(k)}$, $k=1,\dots,m$), we
favor the SRK method DRI1 also for $m>1$.

\begin{table}[tbp]
\renewcommand{\arraystretch}{1.3}
\begin{equation*}
\begin{array}{r|ccccc|ccccc|cccccc}
    0 & & & &&  &&&&&  &&& \\
    \frac{1}{2} & \frac{1}{2} & &&&  & \frac{6-\sqrt{6}}{10} & &&  &&& \\
    1 & -1 & 2 &  &&  & \frac{3+2\sqrt{6}}{5} & & 0 &  &&& \\
    \hline
    0 &&&  &&&&&  &&& \\
    \frac{342}{491} & \frac{342}{491} & & & &  & 3\sqrt{\frac{38}{491}} &&&&  &&& \\
    \frac{342}{491} & \frac{342}{491} & & 0 & &  & -3\sqrt{\frac{38}{491}} && 0 &&  &&& \\
    \hline
    0 & & &&&  &0&&0&&0  &&& \\
    0 & 0 & &&&  & -\frac{214}{513}\sqrt\frac{1105}{991} &&-\frac{491}{513}\sqrt\frac{221}{4955}&&-\frac{491}{513}\sqrt\frac{221}{4955}&&  &&& \\
    0 & 0 & & 0 & & & \frac{214}{513}\sqrt\frac{1105}{991} &&\frac{491}{513}\sqrt\frac{221}{4955} && \frac{491}{513}\sqrt\frac{221}{4955} & & &&& \\
    \hline
    & \frac{1}{6} & & \frac{2}{3} & & \frac{1}{6} &
    \frac{193}{684} & & \frac{491}{1368} & & \frac{491}{1368}
    & & 0 & & \frac16\sqrt{\frac{491}{38}} & & -\frac16\sqrt{\frac{491}{38}}\\
    \cline{2-17}
    & & & & & & -\frac{4955}{7072} & & \frac{4955}{14144}
    & & \frac{4955}{14144} & & 0 & & -\frac18\sqrt{\frac{4955}{221}} & & \frac18\sqrt{\frac{4955}{221}}
\end{array}
\end{equation*}
\caption{Coefficients of the SRK scheme DRI1 with $p_D=3$
and $p_S=2$.} \label{Table-Coeff-DRI1}
\end{table}
\section{Numerical example} \label{Sec:Numerical-Example}
\begin{figure}[tbp]
\begin{center}
\includegraphics[height=5cm]{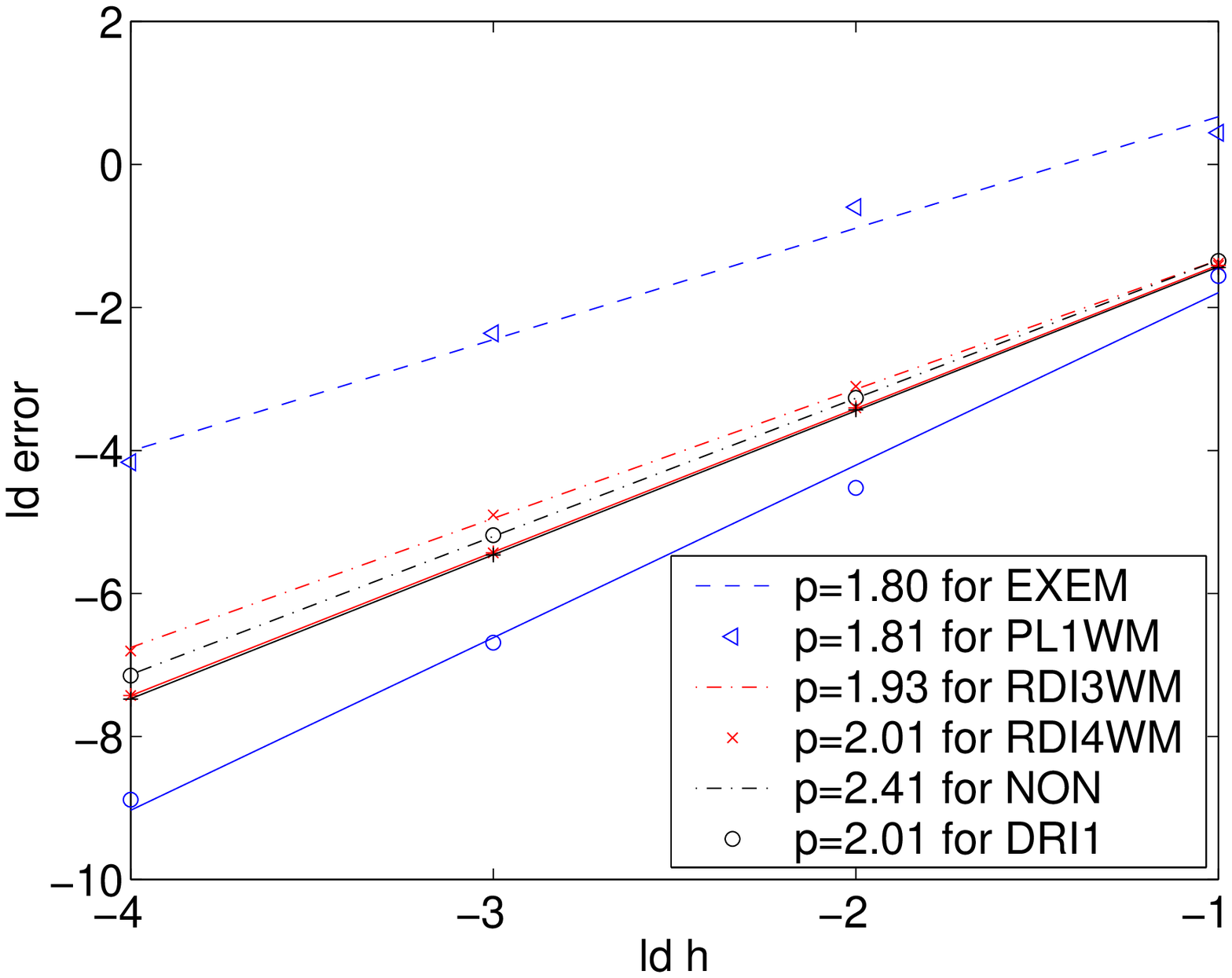}
\includegraphics[height=5cm]{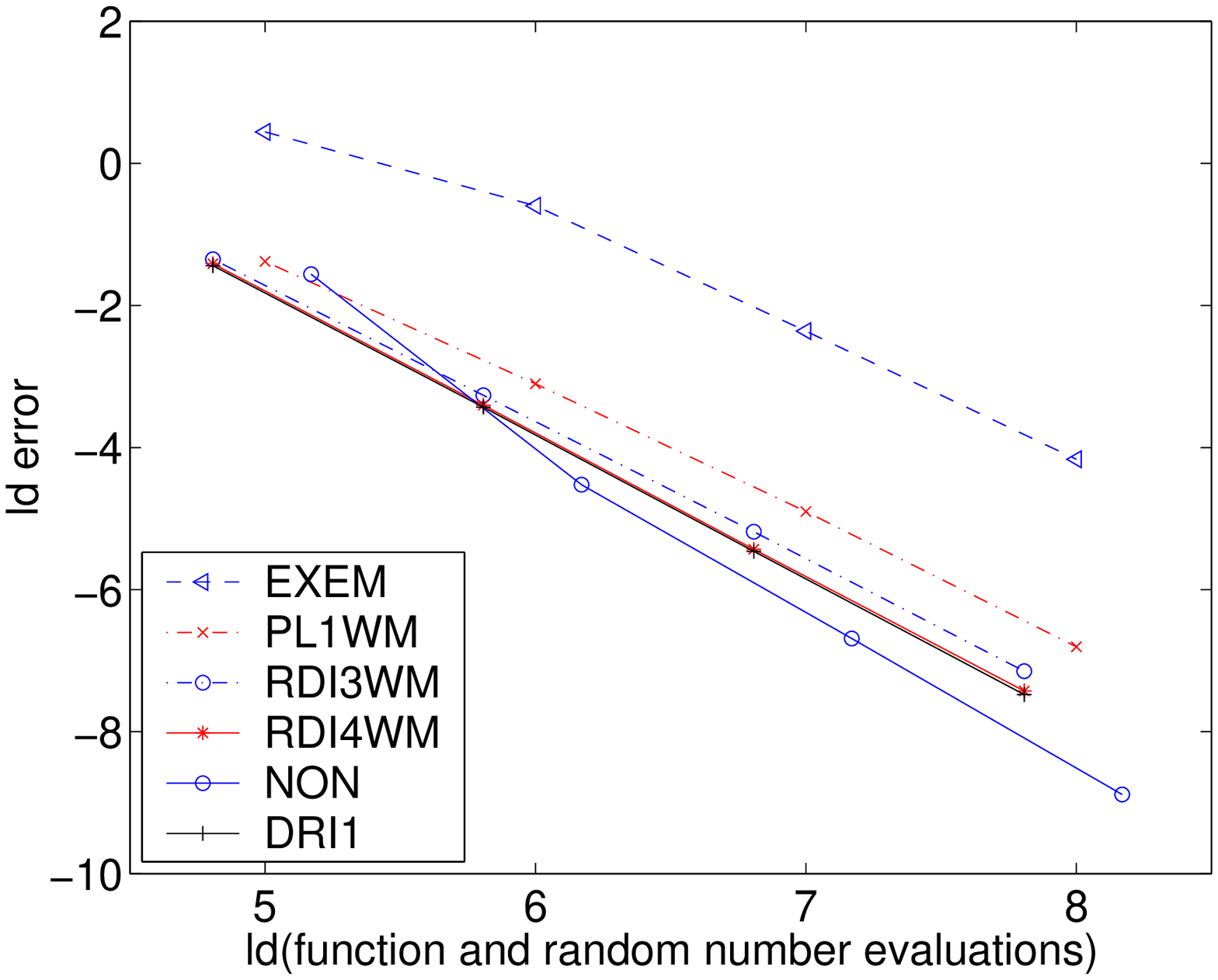}
\caption{Orders of convergence and computational effort per simulation path versus precision for SDE~(\ref{Simu:nonlinear-SDE2}).} \label{Bild001}
\end{center}
\end{figure}

\begin{figure}[tbp]
\begin{center}
\includegraphics[height=5cm]{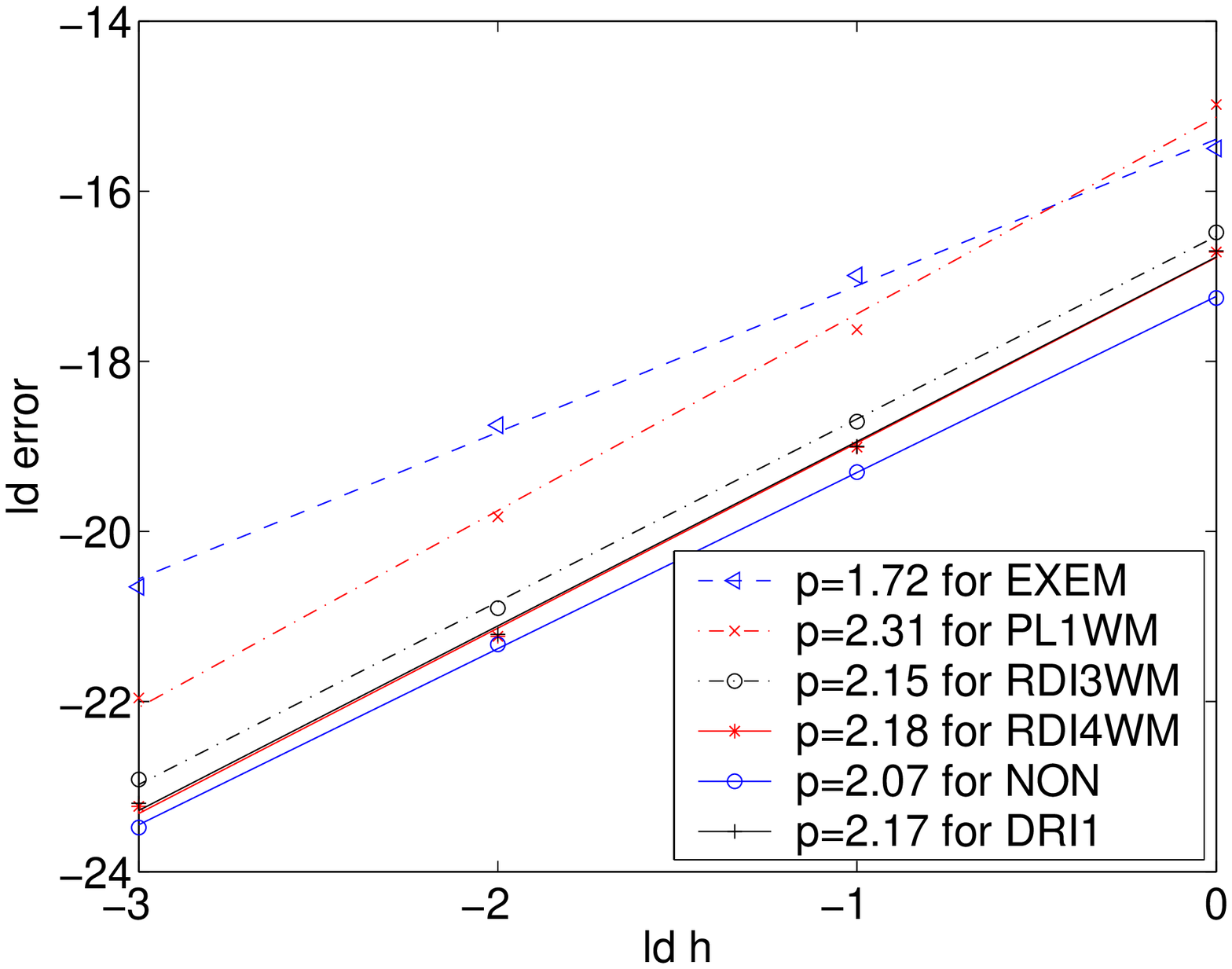}
\includegraphics[height=5cm]{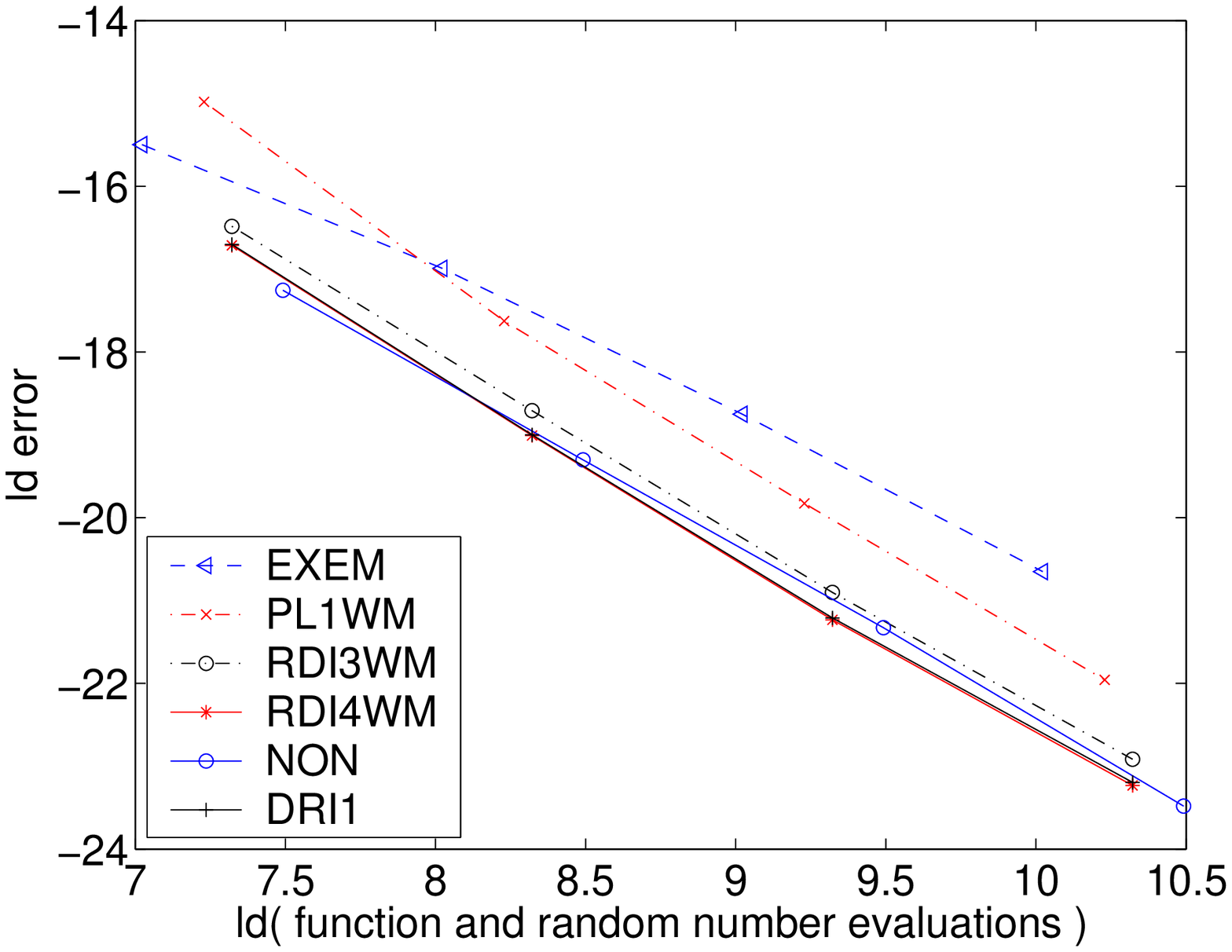}
\caption{Orders of convergence and computational effort per simulation path versus precision for SDE~(\ref{Simu:dm-SDE}).}
\label{Bild002}
\end{center}
\end{figure}

\begin{figure}[tbp]
\begin{center}
\includegraphics[height=5cm]{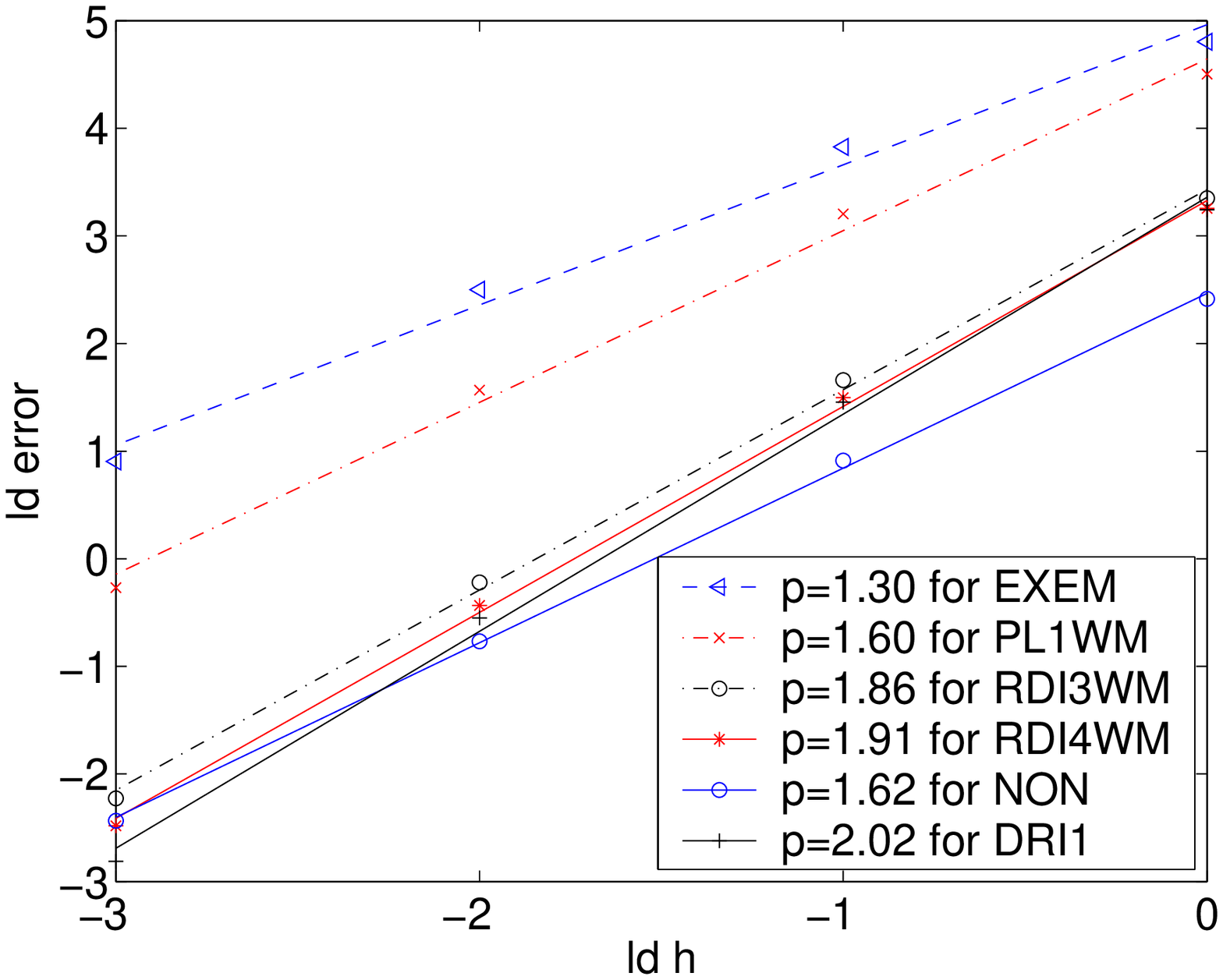}
\includegraphics[height=5cm]{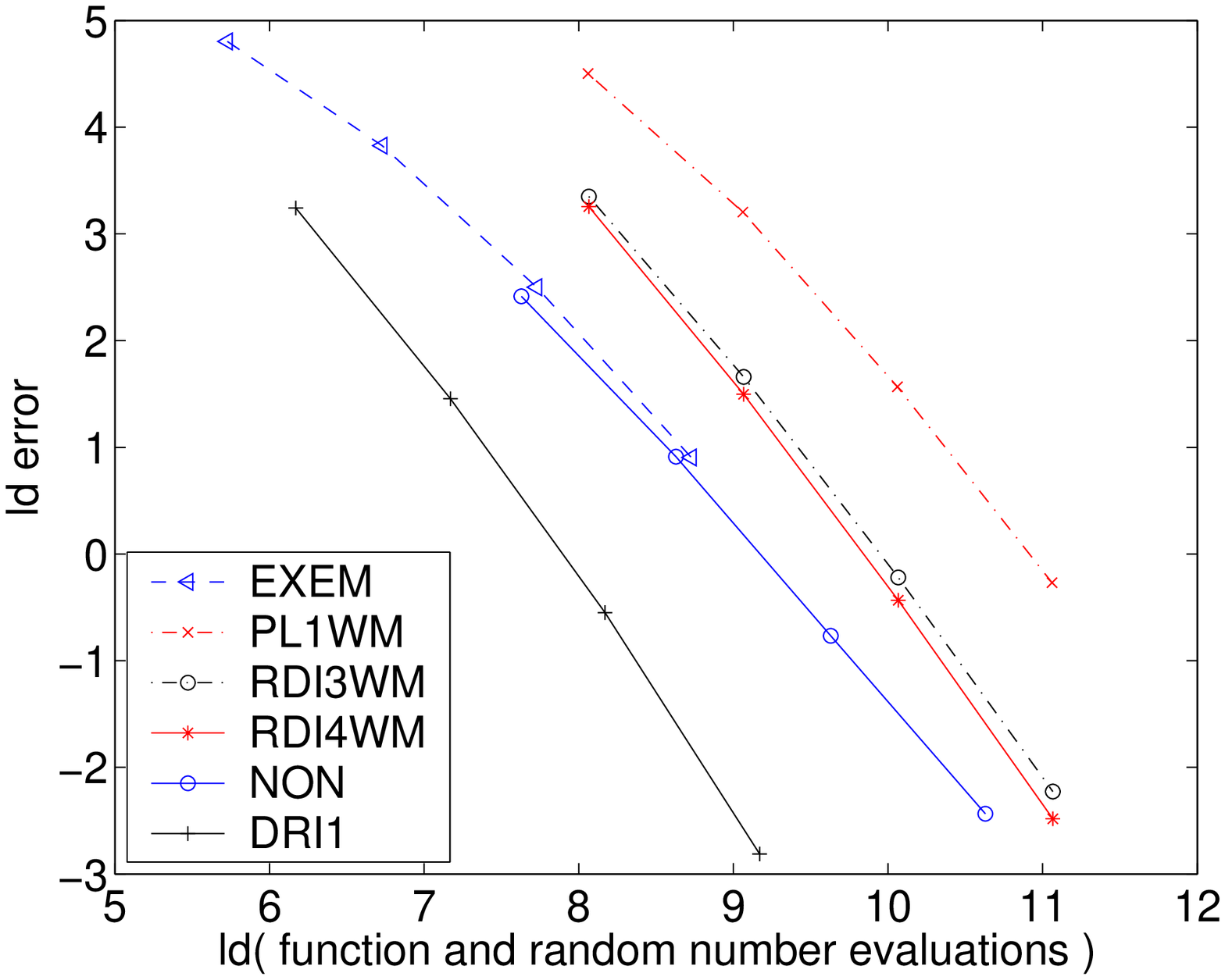}
\caption{Orders of convergence and computational effort per simulation path versus precision for SDE~(\ref{Simu:nonlinear-SDE3}).}
\label{Bild003}
\end{center}
\end{figure}

In the following, the SRK scheme DRI1 presented in
Section~\ref{Sec:Optimal-schemes} is applied to three test equations
in order to analyze its order of convergence in comparison to some
well known schemes.
Therefore, the functional $u = \E(f(X(t)))$ is approximated by a
Monte Carlo simulation. The performance of DRI1 is compared to the
second order SRK schemes PL1WM due to Platen~\cite{KP99}, NON due to Komori \cite{komori07wso}, which in contrast to all other schemes is designed for the weak approximation of Stratonovich SDEs, RDI3WM and RDI4WM due to the authors~\cite{DebrabantRoessler1} and the
ex\-tra\-po\-la\-ted Euler-Maruyama scheme EXEM \cite{TalTub} also
attaining order two, which is given by $2 \E(f(Z^{h/2}(t)))-
\E(f(Z^{h}(t)))$ based on the Euler-Maruyama approximations
$Z^{h/2}(t)$ and $Z^h(t)$ calculated with step sizes $h$ and $h/2$.
The sample average $u_{M,h} = \frac{1}{M} \sum_{k=1}^M
f(Y(t,\omega_k))$, $\omega_k \in \Omega$, of $M$ independent
simulated realizations of the considered approximation $Y(t)$ is
calculated in order to estimate the expectation. In the following,
we denote by $\hat{\mu} = u_{M,h} - \E(f(X(t)))$ the mean error and
by $\hat{\sigma}^2_{\mu}$ the empirical variance of the mean error.
Further, we calculate the confidence interval with boundaries $a$
and $b$ to the level of 90\% for the estimated error $\hat{\mu}$
(see \cite{KP99} for details). \\ \\
As first example, we consider the non-linear
SDE~\cite{KP99,MacNav01}
\begin{equation} \label{Simu:nonlinear-SDE2}
    dX(t) = \left( \tfrac{1}{2} X(t) + \sqrt{X(t)^2 + 1} \right) \,
    dt + \sqrt{X(t)^2 + 1} \, dW(t), \qquad X(0)=0,
\end{equation}
on the time interval $I=[0,2]$ with the solution $X(t) = \sinh (t +
W(t))$. Here, we choose $f(x)=p(\arsinh(x))$, where $p(z) = z^3 -
6z^2 + 8z$ is a polynomial. Then the expectation of the solution
can be calculated as
\begin{equation}
    \E(f(X(t))) = t^3 - 3t^2 + 2t \,\,.
\end{equation}
The solution $\E(f(X(t)))$ is approximated with step sizes $2^{-1},
\ldots, 2^{-4}$ and $M=10^9$ simulations are performed in order to
determine the systematic error of the considered schemes at time
$t=2$.
 The results for the applied schemes are presented in
Table~\ref{Table1}. The orders of convergence correspond to the
slope of the regression lines plotted in the left hand side of
Figure~\ref{Bild001} where we get the order $1.80$ for EXEM, order
$1.81$ for PL1WM, order $1.93$ for RDI3WM, order $2.01$ for RDI4WM, order $2.41$ for NON (applied to the corresponding Stratonovich version of (\ref{Simu:nonlinear-SDE2}))
 and order $2.01$ for the scheme DRI1.\\
Of course, these results have to be related with the computational
effort of the schemes which we take in the following as sum of the
number of evaluations of the drift function $a$ and of each
diffusion function $b^j$, $1 \leq j \leq m$, as well as the number
of random variables that have to be simulated. Then we can compare
the computational effort versus the errors of the analyzed schemes.
The results are presented in the right hand side of
Figure~\ref{Bild001}. The Platen scheme, RDI4WM and the new
scheme DRI1 yield comparable results and all three are better than RDI3WM and much more efficient
than the extrapolated Euler method. For higher precision, NON performs best.
\begin{table}
\caption{Mean errors, empirical variances and confidence intervals
for SDE~(\ref{Simu:nonlinear-SDE2}).} \label{Table1}
\begin{center}
\setlength{\extrarowheight}{-3pt}
\begin{tabular}{c|c|c|c|c|c}
   & $h$ & $\hat{\mu}$ & $\hat{\sigma}_{\mu}^2$  & $a$ & $b$ \\
   \hline
\multirow{4}{*}{EXEM}
   & $2^{-1}$ & -1.359E-00 & 2.990E-06 & -1.359E-00 & -1.359E-00\\
   & $2^{-2}$ & -6.614E-01 & 7.315E-06 & -6.620E-01 & -6.607E-01\\
   & $2^{-3}$ & -1.945E-01 & 8.629E-06 & -1.952E-01 & -1.938E-01\\
   & $2^{-4}$ & -5.570E-02 & 9.014E-06 & -5.641E-02 & -5.499E-02\\
    \hline
\multirow{4}{*}{PL1WM}
   & $2^{-1}$ & -3.837E-01 & 1.885E-06 & -3.841E-01 & -3.834E-01\\
   & $2^{-2}$ & -1.165E-01 & 3.207E-06 & -1.169E-01 & -1.161E-01\\
   & $2^{-3}$ & -3.348E-02 & 2.475E-06 & -3.386E-02 & -3.311E-02\\
   & $2^{-4}$ & -8.949E-03 & 3.447E-06 & -9.390E-03 & -8.509E-03\\
       \hline
\multirow{4}{*}{RDI3WM}
   & $2^{-1}$ & -3.926E-01 & 1.400E-06 & -3.929E-01 & -3.923E-01 \\
   & $2^{-2}$ & -1.041E-01 & 2.787E-06 & -1.045E-01 & -1.037E-01 \\
   & $2^{-3}$ & -2.748E-02 & 2.427E-06 & -2.785E-02 & -2.711E-02 \\
   & $2^{-4}$ & -7.054E-03 & 1.813E-06 & -7.373E-03 & -6.734E-03\\
    \hline
\multirow{4}{*}{RDI4WM}
   & $2^{-1}$ & -3.760E-01 & 1.488E-06 & -3.762E-01 & -3.757E-01 \\
   & $2^{-2}$ & -9.454E-02 & 2.823E-06 & -9.494E-02 & -9.414E-02 \\
   & $2^{-3}$ & -2.318E-02 & 2.441E-06 & -2.355E-02 & -2.281E-02 \\
   & $2^{-4}$ & -5.816E-03 & 1.816E-06 & -6.135E-03 & -5.496E-03\\
       \hline
\multirow{4}{*}{NON}
   & $2^{-1}$ &  -3.393E-01 & 2.530E-06 & -3.396E-01 & -3.389E-01\\
   & $2^{-2}$ &  -4.354E-02 & 3.371E-06 & -4.398E-02 & -4.311E-02\\
   & $2^{-3}$ &  -9.707E-03 & 2.208E-06 & -1.006E-02 & -9.355E-03\\
   & $2^{-4}$ &  -2.119E-03 & 2.952E-06 & -2.526E-03 & -1.711E-03\\
       \hline
\multirow{4}{*}{DRI1}
   & $2^{-1}$ & -3.684E-01 & 1.720E-06 & -3.687E-01 & -3.681E-01\\
   & $2^{-2}$ & -9.271E-02 & 2.939E-06 & -9.312E-02 & -9.231E-02\\
   & $2^{-3}$ & -2.270E-02 & 2.122E-06 & -2.304E-02 & -2.235E-02\\
   & $2^{-4}$ & -5.617E-03 & 2.931E-06 & -6.023E-03 & -5.212E-03
\end{tabular}
\end{center}
\end{table}
%
\\ \\
As a second example, a multi-dimensional SDE with initial value
$X(0)=(1,1)^T$ and noncommutative noise driven by a $2$-dimensional
Wiener process is considered \cite{DebrabantRoessler1}:
\begin{equation} \label{Simu:dm-SDE}
    \begin{split}
    d \begin{pmatrix} X^1 \\ X^2 \end{pmatrix}
    &= \begin{pmatrix} -\frac{273}{512} & 0 \\
    -\frac{1}{160} & -\frac{785}{512}+\frac{\sqrt{2}}{8} \end{pmatrix} \,
    \begin{pmatrix} X^1 \\ X^2 \end{pmatrix} \, dt +
    \begin{pmatrix} \frac{1}{4} X^1 & \frac{1}{16} X^1 \\
    \frac{1-2\sqrt{2}}{4} X^2 & \frac{1}{10} X^1 + \frac{1}{16} X^2
    \end{pmatrix} \, d \begin{pmatrix} W^1 \\ W^2
    \end{pmatrix}.
    \end{split}
\end{equation}
Here, we are interested in the second moments which depend on both,
the drift and the diffusion function (see \cite{KP99} for details).
Therefore, we choose $f(x) = (x^1)^2$
and obtain
\begin{equation}
    \begin{split}
    \E(f(X(t))) = \exp(- t) \, .
    \end{split}
\end{equation}
We approximate $\E(f(X(t)))$ at $t=10$ by $M=8\cdot 10^7$ simulated
trajectories with step sizes $2^{0}, \ldots, 2^{-3}$. The results
for the schemes in consideration are presented in Table~\ref{Table2} and
 Figure~\ref{Bild002}. Here, the order of
convergence is $1.72$ for EXEM, $2.32$ for PL1WM, $2.14$ for RDI3WM, $2.17$ for RDI4WM, $2.07$ for NON and order $2.17$ for our new scheme DRI1.
\begin{table}
\caption{Mean errors, empirical variances and confidence intervals
for SDE~(\ref{Simu:dm-SDE}).} \label{Table2}
\begin{center}
\setlength{\extrarowheight}{-3pt}
\begin{tabular}{c|c|c|c|c|c}
   & $h$ & $\hat{\mu}$ & $\hat{\sigma}_{\mu}^2$  & $a$ & $b$ \\
   \hline
\multirow{4}{*}{EXEM}
   & $2^{-0}$& -2.165E-05 & 1.678E-14 & -2.169E-05 & -2.160E-05\\
   & $2^{-1}$& -7.684E-06 & 1.418E-14 & -7.724E-06 & -7.643E-06\\
   & $2^{-2}$& -2.266E-06 & 1.501E-14 & -2.308E-06 & -2.224E-06\\
   & $2^{-3}$& -6.078E-07 & 3.567E-14 & -6.725E-07 & -5.432E-07\\
    \hline
\multirow{4}{*}{PL1WM}
   & $2^{-0}$ & 3.093E-05 & 9.082E-15 & 3.090E-05 & 3.097E-05\\
   & $2^{-1}$ & 4.947E-06 & 1.085E-14 & 4.906E-06 & 4.987E-06\\
   & $2^{-2}$ & 1.071E-06 & 5.886E-15 & 1.041E-06 & 1.101E-06\\
   & $2^{-3}$ & 2.435E-07 & 4.652E-15 & 2.172E-07 & 2.699E-07\\
    \hline
\multirow{4}{*}{RDI3WM}
   & $2^{-0}$ &-1.092E-05 & 2.481E-15 & -1.094E-05 & -1.090E-05\\
   & $2^{-1}$ &-2.335E-06 & 8.234E-15 & -2.370E-06 & -2.299E-06\\
   & $2^{-2}$ &-5.143E-07 & 5.519E-15 & -5.431E-07 & -4.856E-07\\
   & $2^{-3}$ &-1.285E-07 & 4.581E-15 & -1.546E-07 & -1.023E-07\\
    \hline
\multirow{4}{*}{RDI4WM}
   & $2^{-0}$ &-9.312E-06 & 3.403E-15 & -9.334E-06 & -9.289E-06\\
   & $2^{-1}$ &-1.893E-06 & 8.765E-15 & -1.929E-06 & -1.857E-06\\
   & $2^{-2}$ &-4.096E-07 & 5.591E-15 & -4.386E-07 & -3.807E-07\\
   & $2^{-3}$ &-1.035E-07 & 4.597E-15 & -1.297E-07 & -7.724E-08\\
    \hline
\multirow{4}{*}{NON}
   & $2^{-0}$ &6.396E-06 & 1.588E-14 & 6.347E-06 & 6.445E-06\\
   & $2^{-1}$ &1.548E-06 & 1.266E-14 & 1.504E-06 & 1.591E-06\\
   & $2^{-2}$ &3.799E-07 & 6.172E-15 & 3.495E-07 & 4.102E-07\\
   & $2^{-3}$ &8.544E-08 & 4.713E-15 & 5.889E-08 & 1.120E-07\\
    \hline
\multirow{4}{*}{DRI1}
   & $2^{-0}$ & -9.391E-06 & 3.332E-15 & -9.414E-06 & -9.369E-06\\
   & $2^{-1}$ & -1.908E-06 & 8.710E-15 & -1.944E-06 & -1.872E-06\\
   & $2^{-2}$ & -4.127E-07 & 5.587E-15 & -4.416E-07 & -3.838E-07\\
   & $2^{-3}$ & -1.041E-07 & 4.597E-15 & -1.304E-07 & -7.792E-08
\end{tabular}
\end{center}
\end{table}
%
%
%
%
%
\quad \\ \\
Comparing the computational effort versus precision, in this example
the schemes DRI1 and RDI4WM perform better than NON and RDI3WM and clearly better than EXEM and the Platen scheme.\\
Our last example is a nonlinear SDE with 10 Wiener processes,
\begin{eqnarray}\nonumber
\lefteqn{dX(t)=X(t)\, dt+
\frac1{10}\sqrt{X(t)+\frac12}\, dW_1(t)+
\frac1{15}\sqrt{X(t)+\frac14}\, dW_2(t)}\\
&&+\nonumber
\frac1{20}\sqrt{X(t)+\frac15}\, dW_3(t)+
\frac1{25}\sqrt{X(t)+\frac1{10}}\, dW_4(t)+
\frac1{40}\sqrt{X(t)+\frac1{20}}\, dW_5(t)
\\&&+\nonumber
\frac1{25}\sqrt{X(t)+\frac12}\, dW_6(t)+
\frac1{20}\sqrt{X(t)+\frac14}\, dW_7(t)+
\frac1{15}\sqrt{X(t)+\frac15}\, dW_8(t)\\
&&+\label{Simu:nonlinear-SDE3}
\frac1{20}\sqrt{X(t)+\frac1{10}}\, dW_9(t)+
\frac1{25}\sqrt{X(t)+\frac1{20}}\, dW_{10}(t),\qquad X(0)=1.
\end{eqnarray}
Here, we consider the fourth moment, i.\,e., $f(x)=x^4$, and obtain
\begin{equation}
\E(f(X(t)))=\frac{\scriptstyle 4625768169}{\scriptstyle 73570420483600}-\frac{\scriptstyle 2998776077847}{\scriptstyle 113706563209000}e^{\frac{731453}{360000}t}
+\frac{\scriptstyle 80235120932849}{\scriptstyle 78178246418000}e^{\frac{251453}{60000}t}.
\end{equation}
We approximate $E(f(X(t)))$ at $t=1$ by $M=2\cdot 10^7$ simulated trajectories with step sizes $2^0,\dots,2^{-3}$ and obtain the results presented in Table \ref{Table3} and Figure \ref{Bild003}. The order of convergence is $1.30$ for EXEM, $1.60$ for PL1WM, $1.86$ for RDI3WM, $1.91$ for RDI4WM, $1.62$ for NON and $2.02$ for DRI1. If we take the computational effort into account, we see that DRI1 performs impressively better than all other schemes, which is what we expected for high numbers of Wiener processes.
\begin{table}
\caption{Mean errors, empirical variances and confidence intervals
for SDE~(\ref{Simu:nonlinear-SDE3}).} \label{Table3}
\begin{center}
\setlength{\extrarowheight}{-3pt}
\begin{tabular}{c|c|c|c|c|c}
   & $h$ & $\hat{\mu}$ & $\hat{\sigma}_{\mu}^2$  & $a$ & $b$ \\
   \hline
\multirow{4}{*}{EXEM}
   & $2^{0}$ &-2.793E+01 & 7.005E-04 & -2.794E+01 & -2.792E+01\\
   & $2^{-1}$ &-1.420E+01 & 2.521E-03 & -1.421E+01 & -1.418E+01\\
   & $2^{-2}$ &-5.658E+00 & 7.216E-03 & -5.687E+00 & -5.629E+00\\
   & $2^{-3}$ &-1.872E+00 & 1.040E-02 & -1.907E+00 & -1.837E+00\\
    \hline
\multirow{4}{*}{PL1WM}
   & $2^{-0}$ &-2.266E+01 & 1.183E-03 & -2.268E+01 & -2.265E+01\\
   & $2^{-1}$ &-9.218E+00 & 1.954E-03 & -9.234E+00 & -9.203E+00\\
   & $2^{-2}$ &-2.965E+00 & 4.226E-03 & -2.987E+00 & -2.942E+00\\
   & $2^{-3}$ &-8.294E-01 & 4.294E-03 & -8.519E-01 & -8.070E-01\\
       \hline
\multirow{4}{*}{RDI3WM}
   & $2^{0}$ & -1.019E+01 & 1.727E-03 & -1.021E+01 & -1.018E+01 \\
   & $2^{-1}$ &-3.161E+00 & 2.324E-03 & -3.177E+00 & -3.144E+00 \\
   & $2^{-2}$ &-8.582E-01 & 4.494E-03 & -8.812E-01 & -8.353E-01  \\
   & $2^{-3}$ &-2.136E-01 & 4.373E-03 & -2.363E-01 & -1.910E-01  \\
    \hline
\multirow{4}{*}{RDI4WM}
   & $2^{0}$ & -9.546E+00 & 1.930E-03 & -9.561E+00 & -9.531E+00\\
   & $2^{-1}$ & -2.824E+00 & 2.436E-03 & -2.840E+00 & -2.807E+00\\
   & $2^{-2}$ &-7.398E-01 & 4.557E-03 & -7.629E-01 & -7.167E-01 \\
   & $2^{-3}$ &-1.791E-01 & 4.392E-03 & -2.017E-01 & -1.564E-01 \\
       \hline
\multirow{4}{*}{NON}
   & $2^{0}$  & 5.331E+00 & 4.219E-03 & 5.309E+00 & 5.353E+00\\
   & $2^{-1}$ & 1.883E+00 & 5.097E-03 & 1.858E+00 & 1.907E+00\\
   & $2^{-2}$ & 5.877E-01 & 2.975E-03 & 5.690E-01 & 6.063E-01\\
   & $2^{-3}$ & 1.850E-01 & 2.832E-03 & 1.668E-01 & 2.033E-01\\
       \hline
\multirow{4}{*}{DRI1}
   & $2^{0}$ & -9.465E+00 & 1.103E-03 & -9.476E+00 & -9.453E+00\\
   & $2^{-1}$ &-2.743E+00 & 3.070E-03 & -2.762E+00 & -2.724E+00 \\
   & $2^{-2}$ & -6.834E-01 & 2.531E-03 & -7.006E-01 & -6.662E-01\\
   & $2^{-3}$ & -1.425E-01 & 2.704E-03 & -1.603E-01 & -1.247E-01
\end{tabular}
\end{center}
\end{table}
\section{Conclusion}
\label{Sec:Conclusion}
In the present work, a full classification of the coefficients for a
new class of efficient explicit SRK methods of order $(1,1)$ for
$s=1$ and order $(2,1)$ for $s=2$ stages as well as for order
$(2,2)$ with $s=3$ stages is calculated. Based on this
classification, coefficients for an extension of the deterministic
RK32 scheme to the stochastic case with minimized error constant
are given. For three examples, this scheme is finally compared with
the order two Platen and extrapolated Euler scheme, the schemes RDI3WM and RDI4WM and NON. It turns out
that the new developed scheme performs very well and especially much better than all other schemes in the case of a high number of Wiener processes.
%
%
\section*{Acknowledgements}
The authors are very grateful to the unknown referees for their comments and suggestions.

\end{document}